\documentclass[11pt]{article}
\usepackage{graphicx}
\usepackage{amsmath}
\usepackage{amsfonts}
\usepackage{amssymb}
\newtheorem{theorem}{Theorem}

\newtheorem{proposition}[theorem]{Proposition}

\newenvironment{proof}[1][Proof]{\textbf{#1.} }{\ \rule{0.5em}{0.5em}}

\def\text{\hbox} 

\def\a{\alpha}
\def\b{\beta}

\def\m{\mu}

\def\p{\pi}

\def\s{\sigma}

\def\t{\tau}

\def\f{\phi}

\def\o{\omega}

\def\Si{\Sigma}

\def\F{\Phi}

\def\S{{\bf S}}
\def\E{{\bf E}}
\def\P{\rho}
\def\R{{\cal R}}
\def\A{{\cal A}}
\def\B{{\cal B}}
\def\A'{{\cal A}'}
\def\B'{{\cal B}'}
\def\GG{{\cal G}}
\def\HH{{\cal H}}
\def\FF{{\cal F}}

\def\wh{\widehat}

\def\wti{\widetilde}

\begin{document}

\title{An intrinsic characterization of $p$-symmetric Heegaard splittings
\footnote{{\em An intrinsic characterization of $p$-symmetric
Heegaard splittings}, Proceedings of the Ninth Prague Topological
Symposium, (Prague, 2001), pp.~217--222, Topology Atlas, Toronto,
2002. This contribution is extracted from: M. Mulazzani, {\em On
$p$-symmetric Heegaard splittings}, J. Knot Theory Ramifications
{\bf 9} (2000), no.~8, 1059--1067. Reprinted with permission from
World Scientific Publishing Co.}}
\author{Michele Mulazzani}

\maketitle

\begin{abstract}
{We show that every $p$-fold strictly-cyclic branched covering of
a $b$-bridge link in $\S^3$ admits a $p$-symmetric Heegaard
splitting of genus $g=(b-1)(p-1)$. This gives a complete converse
to a result of Birman and Hilden, and gives an intrinsic
characterization of $p$-symmetric Heegaard splittings as $p$-fold
strictly-cyclic branched coverings of links.
}

\end{abstract}


\section{Introduction}

The concept of $p$-symmetric Heegard splittings has been
introduced by Birman and Hilden (see \cite{BH}) in an extrinsic
way, depending on a particular embedding of the handlebodies of
the splitting in the ambient space $\E^3$. The definition of such
particular splittings was motivated by the aim to prove that every
closed, orientable 3-manifold of Heegaard genus $g\le 2$ is a
2-fold covering of $\S^3$ branched over a link of bridge number
$g+1$ and that, conversely, the 2-fold covering of $\S^3$ branched
over a link of bridge number $b\le 3$ is a closed, orientable
3-manifold of Heegaard genus $b-1$ (compare also \cite{Vi}).

A genus $g$ Heegaard splitting $M=Y_g\cup_{\f}Y'_g$ is called {\it
$p$-symmetric\/}, with $p>1$, if there exist a disjoint embedding
of $Y_g$ and $Y'_g$ into $\E^3$ such that $Y'_g=\t(Y_g)$, for a
translation $\t$ of $\E^3$, and an orientation-preserving
homeomorphism $\P:\E^3\to\E^3$ of period $p$, such that
$\P(Y_g)=Y_g$ and, if $\GG$ denotes the cyclic group of order $p$
generated by $\P$ and $\F:\partial Y_g\to\partial Y_g$ is the
orientation-preserving homeomorphism $\F=\t^{-1}_{\vert\partial
Y'_g}\f$, the following conditions are fulfilled:
\begin{itemize}
\item[i)] $Y_g/\GG$ is homeomorphic to a 3-ball;
\item[ii)] $\mbox{Fix}(\P_{\vert Y_g}^h)=\mbox{Fix}(\P_{\vert Y_g})$,
for each $1\le h\le p-1$;
\item[iii)] $\mbox{Fix}(\P_{\vert Y_g})/\GG$ is an unknotted set of
arcs\footnote{A set of mutually disjoint arcs $\{t_1,\ldots,t_n\}$
properly embedded in a handlebody $Y$ is {\it unknotted\/} if
there is a set of mutually disjoint discs $D=\{D_1,\ldots,D_n\}$
properly embedded in $Y$ such that $t_i\cap D_i=t_i\cap\partial
D_i=t_i$, $t_i\cap D_j=\emptyset$ and $\partial
D_i-t_i\subset\partial Y$ for $1\le i,j\le n$ and $i\neq j$.} in
the ball $Y_g/{\cal G}$;
\item[iv)] there exists an integer $p_0$ such that
$\F\P_{\vert\partial Y_g}\F^{-1}=(\P_{\vert\partial Y_g})^{p_0}$.
\end{itemize}

\medskip

\noindent {\bf Remark 1} By the positive solution of the Smith
Conjecture \cite{MB} it is easy to see that necessarily
$p_0\equiv\pm 1$ mod $p$.

\medskip

The map $\P'=\t\P\t^{-1}$ is obviously an orientation-preserving
homeomorphism of period $p$ of $\E^3$ with the same properties as
$\P$, with respect to $Y'_g$, and the relation
$\f\P_{\vert\partial Y_g}\f^{-1}=(\P'_{\vert\partial Y'_g})^{p_0}$
easily holds.

The {\it $p$-symmetric Heegaard genus\/} $g_p(M)$ of a 3-manifold
$M$ is the smallest integer $g$ such that $M$ admits a
$p$-symmetric Heegaard splitting of genus $g$.

The following results have been established in \cite{BH}:
\begin{enumerate}
\item Every closed, orientable 3-manifold of $p$-symmetric Heegaard
genus $g$ admits a representation as a $p$-fold cyclic covering of
$\S^3$, branched over a link which admits a $b$-bridge
presentation, where $g=(b-1)(p-1)$.
\item The $p$-fold cyclic covering of $\S^3$ branched over a knot of braid number $b$ is a
closed, orientable 3-manifold $M$ which admits a $p$-symmetric
Heegaard splitting of genus $g=(b-1)(p-1)$.
\end{enumerate}

Note that statement 2 is not a complete converse of 1, since it
only concerns knots and, moreover, $b$ denotes the braid number,
which is greater than or equal to (often greater than) the bridge
number. In this paper we fill this gap, giving a complete converse
to statement 1. Since the coverings involved in 1 are
strictly-cyclic (see next section for details on strictly-cyclic
branched coverings of links), our statement will concern this kind
of coverings. More precisely, we shall prove in Theorem
\ref{Theorem 3} that a $p$-fold strictly-cyclic covering of
$\S^3$, branched over a link of bridge number $b$, is a closed,
orientable 3-manifold $M$ which admits a $p$-symmetric Heegaard
splitting of genus $g=(b-1)(p-1)$, and therefore has $p$-symmetric
Heegaard genus $g_p(M)\le (b-1)(p-1)$. This result gives an
intrinsic interpretation of $p$-symmetric Heegaard splittings as
$p$-fold strictly-cyclic branched coverings of links.

\section{Main results}

Let $\b=\{(p_k(t),t)\,\vert\, 1\le k\le
2n\,,\,t\in[0,1]\}\subset\E^2\times[0,1]$ be a geometric
$2n$-string braid of $\E^3$ \cite{Bi}, where
$p_1,\ldots,p_{2n}:[0,1]\to\E^2$ are continuous maps such that
$p_{k}(t)\neq p_{k'}(t)$, for every $k\neq k'$ and $t\in[0,1]$,
and such that
$\{p_1(0),\ldots,p_{2n}(0)\}=\{p_1(1),\ldots,p_{2n}(1)\}$. We set
$P_k=p_k(0)$, for each $k=1,\ldots,2n$, and
$A_i=(P_{2i-1},0),B_i=(P_{2i},0),A'_i=(P_{2i-1},1),B'_i=(P_{2i},1)$,
for each $i=1,\ldots,n$ (see Figure 1). Moreover, we set
$\FF=\{P_1,\ldots,P_{2n}\}$, $\FF_1=\{P_1,P_3\ldots,P_{2n-1}\}$
and $\FF_2=\{P_2,P_4,\ldots,P_{2n}\}$.

The braid $\b$ is realized through an ambient isotopy
${\wh\b}:\E^2\times[0,1]\to\E^2\times[0,1]$,
${\wh\b}(x,t)=(\b_t(x),t)$, where $\b_t$ is an homeomorphism of
$\E^2$ such that $\b_0=\mbox{Id}_{\E^2}$ and $\b_t(P_i)=p_i(t)$,
for every $t\in[0,1]$. Therefore, the braid $\b$ naturally defines
an orientation-preserving homeomorphism
${\wti\b}=\b_1:\E^2\to\E^2$, which fixes the set $\FF$. Note that
$\b$ uniquely defines ${\wti\b}$, up to isotopy of $\E^2$ mod
$\FF$.

Connecting the point $A_i$ with $B_i$ by a circular arc $\a_i$
(called {\it top arc\/}) and the point $A'_i$ with $B'_i$ by a
circular arc $\a'_i$ (called {\it bottom arc\/}), as in Figure 1,
for each $i=1,\ldots,n$, we obtain a $2n$-plat presentation of a
link $L$ in $\E^3$, or equivalently in $\S^3$. As is well known,
every link admits plat presentations and, moreover, a $2n$-plat
presentation corresponds to an $n$-bridge presentation of the
link. So, the bridge number $b(L)$ of a link $L$ is the smallest
positive integer $n$ such that $L$ admits a representation by a
$2n$-plat. For further details on braid, plat and bridge
presentations of links we refer to \cite{Bi}.


\begin{figure}[bht]
 \begin{center}
 \includegraphics*[totalheight=5cm]{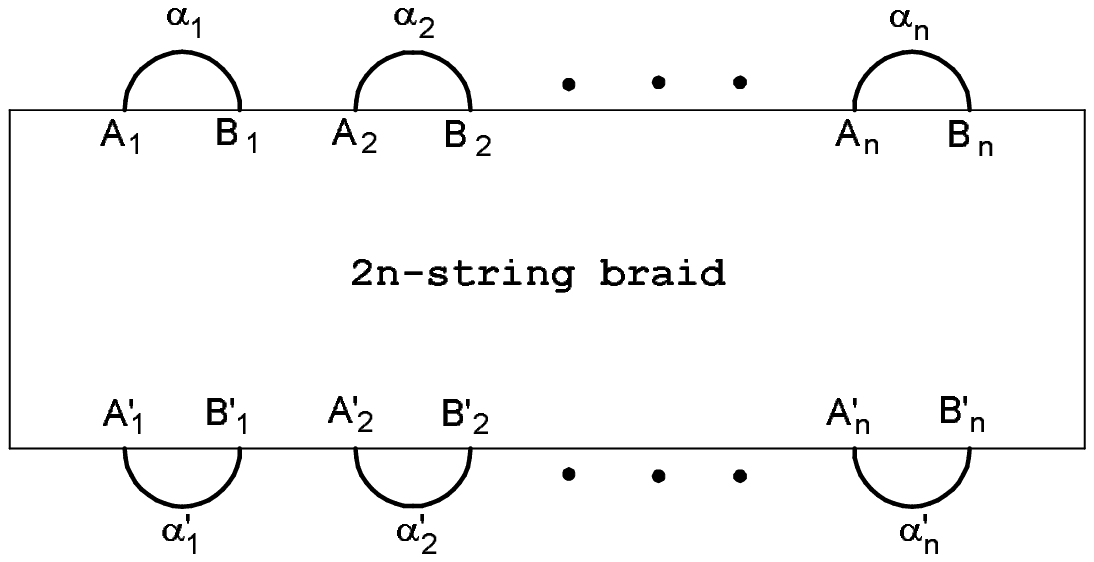}
 \end{center}
 \caption{A $2n$-plat presentation of a link.}

 \label{Fig. 1}

\end{figure}

\medskip

\noindent {\bf Remark 2} A $2n$-plat presentation of a link
$L\subset\E^3\subset\S^3=\E^3\cup\{\infty\}$ furnishes a
$(0,n)$-decomposition \cite{MS}
$(\S^3,L)=(D,A_n)\cup_{\f'}(D',A'_n)$ of the link, where $D$ and
$D'$ are the 3-balls $D=(\E^2\times]-\infty,0])\cup\{\infty\}$ and
$D'=(\E^2\times[1,+\infty[)\cup\{\infty\}$,
$A_n=\a_1\cup\cdots\cup\a_n$, $A'_n=\a'_1\cup\cdots\cup\a'_n$ and
$\f':\partial D\to\partial D'$ is defined by $\f'(\infty)=\infty$
and $\f'(x,0)=({\wti\b}(x),1)$, for each $x\in\E^2$.

\medskip

If a $2n$-plat presentation of a $\m$-component link
$L=\bigcup_{j=1}^{\m}L_j$ is given, each component $L_j$ of $L$
contains $n_j$ top arcs and $n_j$ bottom arcs. Obviously,
$\sum_{j=1}^{\m}n_j=n$. A $2n$-plat presentation of a link $L$
will be called {\it special \/} if:
\begin{itemize}
\item[(1)] the top arcs and the bottom arcs belonging to
$L_1$ are $\a_1,\ldots,\a_{n_1}$ and $\a'_1,\ldots,\a'_{n_1}$
respectively, the top arcs and the bottom arcs belonging to $L_2$
are $\a_{n_1+1},\ldots,\a_{n_1+n_2}$ and
$\a'_{n_1+1},\ldots,\a'_{n_1+n_2}$ respectively, $\ldots$ , the
top arcs and the bottom arcs belonging on $L_\m$ are
$\a_{n_1+\cdots+n_{\m-1}+1},\ldots,\a_{n_1+\cdots+n_{\m}}=\a_{n}$
and
$\a'_{n_1+\cdots+n_{\m-1}+1},\ldots,\a'_{n_1+\cdots+n_{\m}}=\a'_{n}$
respectively;
\item[(2)] $p_{2i-1}(1)\in\FF_1$ and $p_{2i}(1)\in\FF_2$, for each
$i=1,\ldots,n$.
\end{itemize}

It is clear that, because of (2), the homeomorphism $\wti\b$,
associated to a $2n$-string braid $\b$ defining a special plat
presentation, keeps fixed both the sets $\FF_1$ and $\FF_2$.
Although a special plat presentation of a link is a very
particular case, we shall prove that every link admits such kind
of presentation.

\begin{proposition} \label{Proposition special} Every link $L$ admits a special
$2n$-plat presentation, for each $n\ge b(L)$.
\end{proposition}

\begin{proof}
Let $L$ be presented by a $2n$-plat. We show that this
presentation is equivalent to a special one, by using a finite
sequence of moves on the plat presentation which changes neither
the link type nor the number of plats. The moves are of the four
types $I$, $I'$, $II$ and $II'$ depicted in Figure 2. First of
all, it is straightforward that condition (1) can be satisfied by
applying a suitable sequence of moves of type $I$ and $I'$.
Furthermore, condition (2) is equivalent to the following: $(2')$
there exists an orientation of $L$ such that, for each
$i=1,\ldots,n$, the top arc $\a_i$ is oriented from $A_i$ to $B_i$
and the bottom arc $\a'_i$ is oriented from $B'_i$ to $A'_i$.
Therefore, choose any orientation on $L$ and apply moves of type
$II$ (resp. moves of type $II'$) to the top arcs (resp. bottom
arcs) which are oriented from $B_i$ to $A_i$ (resp. from $A'_i$ to
$B'_i$).
\end{proof}


\begin{figure}[bht]
 \begin{center}
 \includegraphics*[totalheight=10cm]{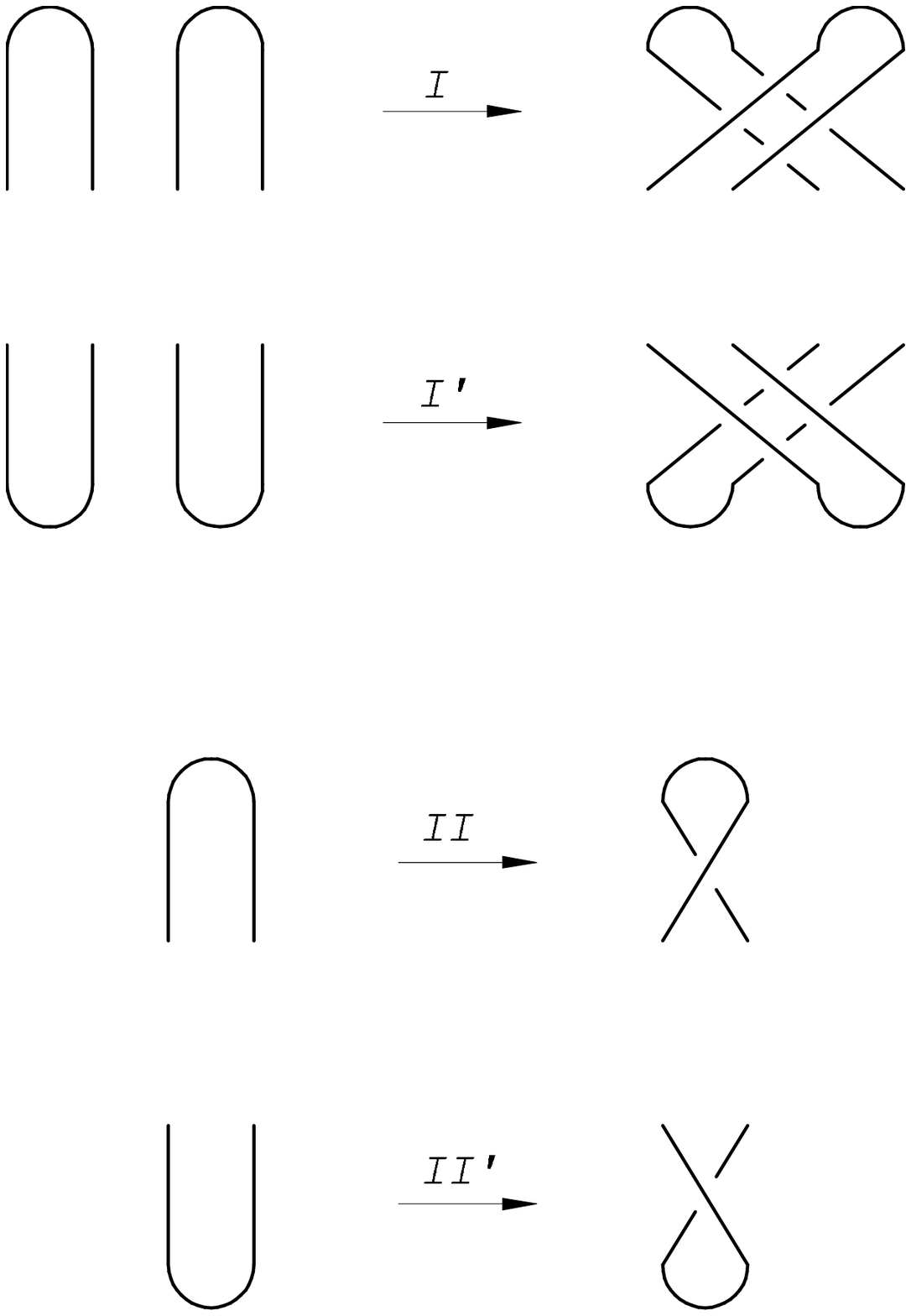}
 \end{center}
 \caption{Moves on plat presentations.}

 \label{Fig. 2}

\end{figure}

\medskip

A $p$-fold branched cyclic covering of an oriented $\m$-component
link $L=\bigcup_{j=1}^{\m}L_j\subset\S^3$ is completely determined
(up to equivalence) by assigning to each component $L_j$ an
integer $c_j\in{\bf Z}_p-\{0\}$, such that the set
$\{c_1,\ldots,c_{\m}\}$ generates the group ${\bf Z}_p$. The
monodromy associated to the covering sends each meridian of $L_j$,
coherently oriented with the chosen orientations of $L$ and
$\S^3$, to the permutation $(1\,2\,\cdots\,p)^{c_j}\in\Si_p$.
Multiplying each $c_j$ by the same invertible element of ${\bf
Z}_p$, we obtain an equivalent covering.

Following \cite{MM} we shall call a branched cyclic covering:
\begin{itemize}
\item [a)] {\it strictly-cyclic\/} if $c_{j'}=c_{j''}$, for every
$j',j''\in\{1,\ldots,\m\}$,
\item [b)] {\it almost-strictly-cyclic\/}
if $c_{j'}=\pm c_{j''}$, for every
$j',j''\in\{1,\ldots,\m\}$,
\item [c)] {\it meridian-cyclic\/} if
$\gcd(b,c_j)=1$, for every $j\in\{1,\ldots,\m\}$,
\item [d)] {\it singly-cyclic\/} if $\gcd(b,c_j)=1$, for some
$j\in\{1,\ldots,\m\}$,
\item [e)] {\it monodromy-cyclic\/} if it is cyclic.
\end{itemize}

The following implications are straightforward: $$\text{ a)
}\Rightarrow\text{ b) }\Rightarrow\text{ c) }\Rightarrow\text{ d)
} \Rightarrow\text{ e) }.$$ Moreover, the five definitions are
equivalent when $L$ is a knot. Similar definitions and properties
also hold for a $p$-fold cyclic covering of a 3-ball, branched
over a set of properly embedded (oriented) arcs.

It is easy to see that, by a suitable reorientation of the link,
an almost-strictly-cyclic covering becomes a strictly-cyclic one.
As a consequence, it follows from Remark 1 that every branched
cyclic covering of a link arising from a $p$-symmetric Heegaard
splitting -- according to Birman-Hilden construction -- is
strictly-cyclic.

Now we show that, conversely, every $p$-fold branched
strictly-cyclic covering of a link admits a $p$-symmetric Heegaard
splitting.

\begin{theorem} \label{Theorem 3} A $p$-fold strictly-cyclic
covering of $\S^3$ branched over a link $L$ of bridge number $b$
is a closed, orientable 3-manifold $M$ which admits a
$p$-symmetric Heegaard splitting of genus $g=(b-1)(p-1)$. So the
$p$-symmetric Heegaard genus of $M$ is
$$g_p(M)\le(b-1)(p-1).$$
\end{theorem}

\begin{proof}
Let $L$ be presented by a special $2b$-plat arising from a braid
$\b$, and let $(\S^3,L)=(D,A_b)\cup_{\f'}(D',A'_b)$ be the
$(0,b)$-decomposition described in Remark 2. Now, all arguments of
the proofs of Theorem 3 of \cite{BH} entirely apply and the
condition of Lemma 4 of \cite{BH} is satisfied, since the
homeomorphism $\wti\b$ associated to $\b$ fixes both the sets
$\FF_1$ and $\FF_2$.
\end{proof}

\medskip

As a consequence of Theorem \ref{Theorem 3} and Birman-Hilden
results, there is a natural one-to-one correspondence between
$p$-symmetric Heegaard splittings and $p$-fold strictly-cyclic
branched coverings of links.

\vspace{15 pt} {MICHELE MULAZZANI, Department
of Mathematics, University of Bologna, I-40127 Bologna, ITALY,
and C.I.R.A.M., Bologna, ITALY. E-mail: mulazza@dm.unibo.it}

\end{document}